\documentclass[12pt]{article}
\usepackage[final]{epsfig}
\usepackage{graphics}
\usepackage{amsmath}
\usepackage{amsfonts}
\usepackage{latexsym}
\usepackage{amssymb}
\usepackage{graphicx}
\usepackage{epstopdf}
\def\rank{\mathop{\rm rank}}
\def\card{\mathop{\rm card}}
\def\Ind{\mathop{\rm Ind}}
\def\Res{\mathop{\rm Res}}
\def\Lie{\mathop{\rm Lie}}
\input pictexwd.tex
\font\tiny=cmr5
\font\small=cmr10

\newtheorem{proposition}{Proposition}[section]
\newtheorem{theorem}[proposition]{Theorem}

\newtheorem{corollary}[proposition]{Corollary}

\title{Jordan types of triangular matrices over a finite field}

\author{Dmitry Fuchs, Alexandre Kirillov sr.}

\date{\it To the memory of Louiza Kirillova (1937-2021)}

\begin{document}

\maketitle

\begin{abstract}
{\small Let $\lambda$ be a partition of an integer $n$ and ${\mathbb F}_q$ be a finite field of order $q$. Let $P_\lambda(q)$ be the number of strictly upper triangular $n\times n$ matrices of the Jordan type $\lambda$. It is known that the polynomial $P_\lambda$ has a tendency to be divisible by high powers of $q$ and $Q=q-1$, and we put $P_\lambda(q)=q^{d(\lambda)}Q^{e(\lambda)}R_\lambda(q)$, where $R_\lambda(0)\neq0$ and $R_\lambda(1)\neq0$. In this article, we study the polynomials $P_\lambda(q)$ and $R_\lambda(q)$. Our main results: an explicit formula for $d(\lambda)$ (an explicit formula for $e(\lambda)$ is known, see Proposition 3.3 below), a recursive formula for $R_\lambda(q)$ (a similar formula for $P_\lambda(q)$ is known, see Proposition 3.1 below), the stabilization of $R_\lambda$ with respect to extending $\lambda$ by adding strings of 1's, and an explicit formula for the limit series $R_{\lambda1^\infty}$. Our studies are motivated by projected applications to the orbit method in the representation theory of nilpotent algebraic groups over finite fields.}
\end{abstract}

\section{Introduction}\label{intro}

This paper is a part of a bigger program of the application of the orbit method to the representation theory of triangular matrix group. There are several  reasons for this endeavour. 

It is known (see, e.g. \cite{K1}) that for a wide class of Lie groups $G$ the explicit and transparent answer to the main questions of the representation theory may be formulated, and in some cases proved, in terms of ``coadjoint orbits," that is, orbits of the natural action of $G$ in the dual space $\mathfrak g^\ast$ of the Lie algebra ${\mathfrak g}=\Lie G$.

There is a hope that some modification of the orbit method will work also for the group $N_n({\mathbb F}_q)$ of unitriangular $n\times n$ matrices over the finite field ${\mathbb F}_q$. A simple argument shows that the number of coadjoint orbits of the group $\mathfrak g^\ast$ is equal to the number of adjoint orbits of this group in the Lie algebra ${\mathfrak n}_+(n)$ of strictly upper triangular $n\times n$ matrices over ${\mathbb F}_q$. (Namely, the Fourier transform establishes an isomorphism between the space of functions on ${\mathfrak n}_+(n)$ constant on adjoint orbits and the space of functions on ${\mathfrak n}_+(n)^\ast$ constant on coadjoint orbits.) But this number is the same as the number conjugacy classes of the group $N_n({\mathbb F}_q)$, that is the number of irreducible representations of this group. However, no explicit constriction of a representation, corresponding to a coadjoint orbit, is known, and even the sufficient information of the number of these orbits does not exist. 

Along with the partition of the algebra ${\mathfrak n}_+(n)$ into the adjoint orbits, there exists a more rough partition into Jordan types, which correspond to partitions of $n$; we denote the set of these partition by ${\mathcal P}(n)$. For a partition $\lambda=(\lambda_1,\lambda_2,\ldots,\lambda_N), \lambda_1\ge\lambda_2\ge\ldots\ge\lambda_N$ of $n$ we use also the ``symbolic notation" $1^{\alpha_1}2^{\alpha_2}\dots N^{\alpha_N}$ where $\alpha_j=\card\{i\mid \lambda_i=j\}$; usually (but not always) factors with $\alpha_j=0$ as well as the exponent $\alpha_j=1$, are omitted.

Notice that for a matrix $A\in{\mathfrak n}_+(n)$ the sequence $\alpha_1,\alpha_2,\alpha_3,\dots$ of exponents in the symbolic notation of its Jordan type $\lambda=1^{\alpha_1}2^{\alpha_2}3^{\alpha_3}\dots$ (with zeroes not removed) can be calculated as the sequence of second differences of the rank sequence $\{r_k(A)=\rank A^k\mid k=0,1,2,\dots$, that is $$\alpha_k=\frac{r_{k-1}(A)-2r_k(A)+r_{k+1}(A)}2.$$

Let $P_\lambda(q)$ be the number of matrices  from ${\mathfrak n}_n({\mathbb F}_q)$ of Jordan type $\lambda\in{\mathcal P}(n)$. Thus, $$\sum_{\lambda\in{\mathcal P}(n)}P_\lambda(q)=q^{\frac{n(n-1)}2}.$$It is well known (and explained in Section 3 below) that $P_\lambda$ is a polynomial of the variable $q$. This polynomials were studied in many works (see, for example, ..... ). It is known, in particular that these polynomials have a tendency to be divisible by high degrees of $q$ and $Q=q-1$. There arises a polynomial $R_\lambda$ defined by the formula$$P_\lambda(q)=q^{d(\lambda)}Q^{e(\lambda)}R_\lambda(q),\ {\rm where}\ R_\lambda(0)\ne0,R_\lambda(1)\ne0.$$. This article is our contribution to the study of polynomials $P_\lambda$ and $R_\lambda$. Our main results: an explicit formula for $d(\lambda)$; stabilization for $k\to\infty$ of polynomials $R_{\lambda1^k}$ where $\lambda1^k$ is a partition of $n+k$ obtained from $\lambda\in{\mathcal P}(n)$ by adding $k$ ones; an explicit formula for the limit series $R_{\lambda1^\infty}$. 

The plan of the article is the following. Section 2 contains the traditional material concerning partitions and Young diagrams; our main goal there is to establish notations used in the rest of the article. Section 3 contains results, which we consider as known, namely, a recursive formula for $P_\lambda$ and explicit formulas for $\deg P_\lambda$ and for $e(\lambda)$. For the last two, we prefer to provide short proofs. Section 4 is devoted to computing $d(\lambda)$. In Section 5 we derive some immediate corollaries for $R_\lambda$. Section 6 is devoted to the stabilization phenomenon for $R_\lambda$. In Appendix, we provide a table for polynomials $P_\lambda$ and $R_\lambda$ for all $\lambda\in{\mathcal P}(n),\ n\le10$; in a slightly different form, this table was compiled by D. Golubenko.\smallskip

{\bf Acknowledgements.} We are grateful to D. Golubenko for his help with a table for $P_\lambda$.

\section{Partitions and Young diagrams.}\label{partitions}

A {\it partition} $\lambda$ of a positive integer $n$ is a sequence $\lambda_1,\lambda_2,\ldots,\lambda_N$ of integers such that $\lambda_1\ge\lambda_2\ge\ldots\ge\lambda_N>0$ and $\lambda_1+\lambda_2+\dots+\lambda_N=n$. The set of all partition of $n$ is denoted as ${\cal P}()n)$. Graphically, a partition $\lambda=(\lambda_1,\lambda_2,\dots,\lambda_N)$ of $n$ is represented by {\it Young diagrams} like the one shown in Figure 1 below.

\centerline{
\beginpicture
\setcoordinatesystem units <1in,1in> point at 0 0
\linethickness=.6pt
\setplotsymbol({\small.})
\plot 0 -1.8 0 0 1 0 /
\plot 0 -.2 1 -.2 /
\plot 0 -.4 1 -.4 1 0 /
\plot 0 -.6 .6 -.6 /
\plot 0 -.8 .6 -.8 /
\plot 0 -1 .6 -1 .6 0 /
\plot 0 -1.2 .4 -1.2 /
\plot 0 -1.4 .4 -1.4 .4 0 /
\plot 0 -1.6 .2 -1.6 /
\plot 0 -1.8 .2 -1.8 .2 0 /
\plot 0 -1.8 .2 -1.8 .2 0 /
\plot .8 0 .8 -.4 /
\put{$\lambda_1$} at -.12 -.1
\put{$\lambda_2$} at -.12 -.3
\put{$\lambda_3$} at -.12 -.5
\put{$\lambda_4$} at -.12 -.7
\put{$\lambda_5$} at -.12 -.9
\put{$\lambda_6$} at -.12 -1.1
\put{$\lambda_7$} at -.12 -1.3
\put{$\lambda_8$} at -.12 -1.5
\put{$\lambda_9$} at -.12 -1.7
\put{$\lambda_N=$} at -.42 -1.7
\put{$\lambda'_1$} at .1 .12
\put{$\lambda'_2$} at .3 .12
\put{$\lambda'_3$} at .5 .12
\put{$\lambda'_4$} at .7 .12
\put{$\lambda'_5$} at .9 .12
\put{$=\lambda'_{N'}$} at 1.2 .12
\put{$(i_1,i'_1)$} at 1.27 -.3
\put{$(i_2,i'_2)$} at .87 -.9
\put{$(i_3,i'_3)$} at .67 -1.3
\put{$(i_4,i'_4)$} at .47 -1.7
\setplotsymbol({\tiny.})
\plot .84 -.2 .8 -.24 /
\plot .88 -.2 .8 -.28 /
\plot .92 -.2 .8 -.32 /
\plot .96 -.2 .8 -.36 /
\plot 1 -.2 .8 -.4 /
\plot 1 -.24 .84 -.4 /
\plot 1 -.28 .88 -.4 /
\plot 1 -.32 .92 -.4 /
\plot 1 -.36 .96 -.4 /
\plot .44 -.8 .4 -.84 /
\plot .48 -.8 .4 -.88 /
\plot .52 -.8 .4 -.92 /
\plot .56 -.8 .4 -.96 /
\plot .6 -.8 .4 -1 /
\plot .6 -.84 .44 -1 /
\plot .6 -.88 .48 -1 /
\plot .6 -.92 .52 -1 /
\plot .6 -.96 .56 -1 /
\plot .24 -1.2 .2 -1.24 /
\plot .28 -1.2 .2 -1.28 /
\plot .32 -1.2 .2 -1.32 /
\plot .36 -1.2 .2 -1.36 /
\plot .4 -1.2 .2 -1.4 /
\plot .4 -1.24 .24 -1.4 /
\plot .4 -1.28 .28 -1.4 /
\plot .4 -1.32 .32 -1.4 /
\plot .4 -1.36 .36 -1.4 /
\plot .04 -1.6 0 -1.64 /
\plot .08 -1.6 0 -1.68 /
\plot .12 -1.6 0 -1.72 /
\plot .16 -1.6 0 -1.76 /
\plot .2 -1.6 0 -1.8 /
\plot .2 -1.64 .04 -1.8 /
\plot .2 -1.68 .08 -1.8 /
\plot .2 -1.72 .12 -1.8 /
\plot .2 -1.76 .16 -1.8 /
\endpicture
\beginpicture
\setcoordinatesystem units <1in,1in> point at 0 0
\linethickness=.6pt
\setplotsymbol({\small.})
\plot 0 -1.8 0 0 1 0 /
\plot 0 -.2 1 -.2 /
\plot 0 -.4 1 -.4 1 0 /
\plot 0 -.6 .6 -.6 /
\plot 0 -.8 .6 -.8 /
\plot 0 -1 .6 -1 .6 0 /
\plot 0 -1.2 .4 -1.2 /
\plot 0 -1.4 .4 -1.4 .4 0 /
\plot 0 -1.6 .2 -1.6 /
\plot 0 -1.8 .2 -1.8 .2 0 /
\plot 0 -1.8 .2 -1.8 .2 0 /
\plot .8 0 .8 -.4 /
\put{$\left\{\displaystyle\frac\strut\strut\right.$} at -.06 -.2
\put{$\left\{\begin{array}{c}\strut\\ \strut\\ \strut\end{array}\right.$} at -.01 -.7
\put{$\left\{\displaystyle\frac\strut\strut\right.$} at -.06 -1.2
\put{$\left\{\displaystyle\frac\strut\strut\right.$} at -.06 -1.6
\put{$\alpha_5$} at -.24 -.2
\put{$\alpha_3$} at -.24 -.7
\put{$\alpha_2$} at -.24 -1.2
\put{$\alpha_1$} at -.24 -1.6
\setplotsymbol({\bf.})
\plot 0 0 1 0 /
\plot 0 -.4 1 -.4 /
\plot 0 -1 .6 -1 /
\plot 0 -1.4 .4 -1.4 /
\plot 0 -1.8 .2 -1.8 /
\setplotsymbol({\tiny.})
\plot .84 -.2 .8 -.24 /
\plot .88 -.2 .8 -.28 /
\plot .92 -.2 .8 -.32 /
\plot .96 -.2 .8 -.36 /
\plot 1 -.2 .8 -.4 /
\plot 1 -.24 .84 -.4 /
\plot 1 -.28 .88 -.4 /
\plot 1 -.32 .92 -.4 /
\plot 1 -.36 .96 -.4 /
\plot .44 -.8 .4 -.84 /
\plot .48 -.8 .4 -.88 /
\plot .52 -.8 .4 -.92 /
\plot .56 -.8 .4 -.96 /
\plot .6 -.8 .4 -1 /
\plot .6 -.84 .44 -1 /
\plot .6 -.88 .48 -1 /
\plot .6 -.92 .52 -1 /
\plot .6 -.96 .56 -1 /
\plot .24 -1.2 .2 -1.24 /
\plot .28 -1.2 .2 -1.28 /
\plot .32 -1.2 .2 -1.32 /
\plot .36 -1.2 .2 -1.36 /
\plot .4 -1.2 .2 -1.4 /
\plot .4 -1.24 .24 -1.4 /
\plot .4 -1.28 .28 -1.4 /
\plot .4 -1.32 .32 -1.4 /
\plot .4 -1.36 .36 -1.4 /
\plot .04 -1.6 0 -1.64 /
\plot .08 -1.6 0 -1.68 /
\plot .12 -1.6 0 -1.72 /
\plot .16 -1.6 0 -1.76 /
\plot .2 -1.6 0 -1.8 /
\plot .2 -1.64 .04 -1.8 /
\plot .2 -1.68 .08 -1.8 /
\plot .2 -1.72 .12 -1.8 /
\plot .2 -1.76 .16 -1.8 /
\endpicture
\hskip1.1in
\beginpicture
\setcoordinatesystem units <1in,1in> point at 0 0
\linethickness=.6pt
\setplotsymbol({\small.})
\plot 0 -1.8 0 0 1 0 /
\plot 0 -.2 1 -.2 /
\plot 0 -.4 1 -.4 1 0 /
\plot 0 -.6 .6 -.6 /
\plot 0 -.8 .6 -.8 /
\plot 0 -1 .6 -1 .6 0 /
\plot 0 -1.2 .4 -1.2 /
\plot 0 -1.4 .4 -1.4 .4 0 /
\plot 0 -1.6 .2 -1.6 /
\plot 0 -1.8 .2 -1.8 .2 0 /
\plot 0 -1.8 .2 -1.8 .2 0 /
\plot .8 0 .8 -.4 /
\put{$\overbrace{\hskip.4in}$} at .8 .06
\put{$\alpha'_2$} at .8 .2
\put{$\alpha'_5$} at .5 .12
\put{$\alpha'_7$} at .3 .12
\put{$\alpha'_9$} at .1 .12
\setplotsymbol({\bf.})
\plot 0 0 0 -1.8 /
\plot .2 0 .2 -1.8 /
\plot .4 0 .4 -1.4 /
\plot .6 0 .6 -1 /
\plot 1 0 1 -.4 /
\setplotsymbol({\tiny.})
\plot .84 -.2 .8 -.24 /
\plot .88 -.2 .8 -.28 /
\plot .92 -.2 .8 -.32 /
\plot .96 -.2 .8 -.36 /
\plot 1 -.2 .8 -.4 /
\plot 1 -.24 .84 -.4 /
\plot 1 -.28 .88 -.4 /
\plot 1 -.32 .92 -.4 /
\plot 1 -.36 .96 -.4 /
\plot .44 -.8 .4 -.84 /
\plot .48 -.8 .4 -.88 /
\plot .52 -.8 .4 -.92 /
\plot .56 -.8 .4 -.96 /
\plot .6 -.8 .4 -1 /
\plot .6 -.84 .44 -1 /
\plot .6 -.88 .48 -1 /
\plot .6 -.92 .52 -1 /
\plot .6 -.96 .56 -1 /
\plot .24 -1.2 .2 -1.24 /
\plot .28 -1.2 .2 -1.28 /
\plot .32 -1.2 .2 -1.32 /
\plot .36 -1.2 .2 -1.36 /
\plot .4 -1.2 .2 -1.4 /
\plot .4 -1.24 .24 -1.4 /
\plot .4 -1.28 .28 -1.4 /
\plot .4 -1.32 .32 -1.4 /
\plot .4 -1.36 .36 -1.4 /
\plot .04 -1.6 0 -1.64 /
\plot .08 -1.6 0 -1.68 /
\plot .12 -1.6 0 -1.72 /
\plot .16 -1.6 0 -1.76 /
\plot .2 -1.6 0 -1.8 /
\plot .2 -1.64 .04 -1.8 /
\plot .2 -1.68 .08 -1.8 /
\plot .2 -1.72 .12 -1.8 /
\plot .2 -1.76 .16 -1.8 /
\endpicture
}\vskip.2in

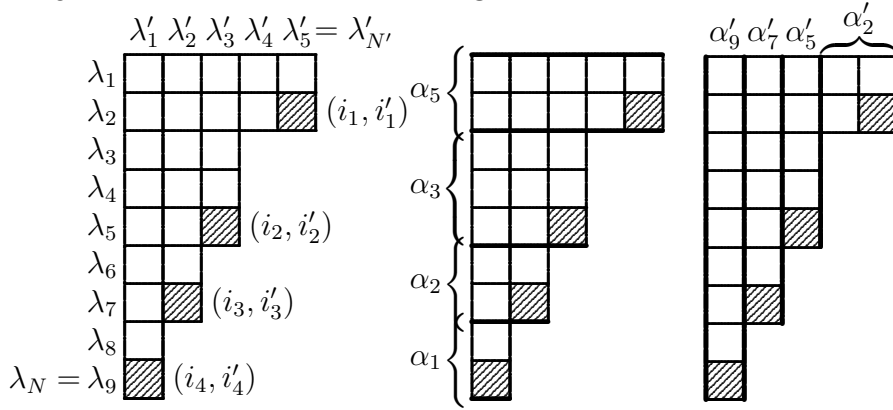
\centerline{Figure 1. Young diagram of the partition (5,5,3,3,3,2,2,1,1) of 25.}

\vskip.2in

The number $n$ (25 in Figure 1) is the number of square cells of the diagram. The ``parts'' $\lambda_1,\lambda_2,\dots,\lambda_N$ are the lengths of rows (counted from the top to the bottom). Associated with the partition $\lambda$, there is the {\it dual} partition $\lambda'$ of the same $n$. Its parts $\lambda'_1,\lambda'_2,\ldots,\lambda'_{N'}$ are the heights of the columns of the same diagram. In the language of formulas, it can be said that $\lambda'_j=\card\{i\mid \lambda_i\ge j\}$. Also, we can say that the Young diagram of the partition $\lambda'$ can be obtained from the Young diagram of $\lambda$ by the reflection in the bisector.

Notice that $N=\lambda'_1$ and $N'=\lambda_1$.

Every cell of a Young diagram has ``coordinates'' $(i,i')$, which are, respectively the numbers of the row and the column to which this cell belongs. 

A cell of a Young diagram is called {\it removable}, if after its removing the Young diagram remains a Young diagram. Geometrically, removable cells correspond to the outer angles of the polygonal line, which bounds the diagram from the right (or from the bottom). The Young diagram in Figure 1 has 4 removable cells; they are shadowed in the picture. If there are $s$ removable cells, then we denote their coordinates as $(i_1,i'_1),(i_2,i'_2),\dots,(i_s,i'_s)$ in such a way that $i_1<i_2<\ldots<i_s\, (=N)$ and $(N'=)\, i'_1>i'_2>\ldots>i'_s\, (=1)$. In Figure 1, $s=4$ and coordinates of the removable cells are $(2,5), (5,3), (7,2), (9,1)$. \smallskip

The partition $\lambda$ in Figure 1 may be described as $1^22^23^35^2$, while the partition $\lambda'$ is described as $2^25^17^19^1$ or $2^2579$. Certainly, it is also possible to order the terms in the sequence $i_1^{\alpha_1}i_2^{\alpha_2}\dots i_s^{\alpha_s}$ in the order of descending $i_k$. For example, there are 11 partitions of $n=6$: $$\begin{array} {c} 6;\ 51;\ 42;\ 41^2;\ 3^2;\ 321;\ 31^3;\ 2^3;\ 2^21^2; 21^4; 1^6;\\ 6'=1^6; (51)'=21^4; (42)'=2^21^2; (41^2)'=31^3; (3^2)'=2^3;\ (321)'=321\end{array}$$

The {\it non-zero} numbers $\alpha_k$, as well as the similar non-zero numbers $\alpha'_k$ for the dual partition $\lambda'$, also have a clear geometric sense. Namely, the Young diagram of the patition may be subdivided into {\it horizontal strips} composed of rows of equal lengths and into {\it vertical strips} composed of columns of equal heights. Both these subdivision are shown in Figure 1. (Notice that every horizontal strip, as well as every vertical strip, contains precisely one removable cell, which shows that the number of horizontal strips, as well as the number of vertical strips, coincides with the number of removable cells, so these two numbers are equal to each other.) And it is clear that $\alpha_k$ is the height of the horizontal strips composed of rows of the length $k$, while $\alpha'_k$ is the width of the vertical strip composed of columns of the height $k$. By the way, this creates a simple way to deduce the last description of the partitions $\lambda$ and $\lambda'$ from the coordinates of removable cells. Namely, $$\begin{array} {c} \lambda=(i'_1)^{i_s-i_{s-1}}(i'_2)^{i_{s-1}-i_{s-2}}\ldots(i'_{s-1})^{i_2-i_1}(i'_s)^{i_1},\\ \lambda'=(i_1)^{i'_1-i'_2}(i_2)^{i'_2-i'_3}\ldots(i_{s-1})^{i'_{s-1}-i'_s}(i_s)^{i'_s-1};\end{array}\eqno(1)$$ in particular, such presentation of partitions dual to each other have equal lengths (the number of ``factors'').

There are two more formulas, which are worth mentioning (in these formulas we {\it do not omit} zero $\alpha$ and $\alpha'$): $$\lambda_k=\alpha'_k+\alpha'_{k+1}+\alpha'_{k+2}+\dots\ {\rm and}\ \lambda'_k=\alpha_k+\alpha_{k+1}+\alpha_{k+2}+\dots.\eqno(2)$$

In conclusion, remind that partition and Young diagrams play very special role in the theory of representations of symmetric groups. First of all, in the symmetric group ${\mathfrak S}_n$, the conjugacy classes are labelled by  partitions of $n$: a conjugacy class corresponds to the partition of $n$ into the sum of lengths of cycles of permutations from this class. Furthermore, according to the classical theorem of Frobenius, the number of (isomorphism classes of) irreducible representations of a finite group is equal to the number of conjugacy classes in this group. However, in the case of the symmetric group ${\mathfrak S}_n$, there exists a canonical bijection between the set of conjugacy classes (which is the same as ${\mathcal P}(n)$) and the set of classes of irreducible representations. Thus, the irreducible representations of ${\mathfrak S}_n$ may be labelled as $V_\lambda$ where $\lambda$ is a partition of $n$. There are multiple explicit construction of $V_\lambda$, as well as multiple explicit formulas for $\dim V_\lambda$ [references]. We will mention here a theorem from the representation theory, which we will need below. Let the Young diagram of the partition $\lambda$ of $n$ has $s$ removable cells, and let $\lambda\hskip-3pt\downarrow_j\ (1\le j\le s)$ be partition of $n-1$ whose Young diagrams are obtained from the Young diagram of $\lambda$ by removing the removable cell $(i_j,i'_j)$. (Notice that the operation of removing the cell $(i_j,i'_j)$ makes the following change in the sequence $\{\alpha_k\}$: $\alpha_j\ne0$ looses 1, and $\alpha_{j-1}$, which could be zero, gains 1.) Then $V_\lambda$ regarded as a representation of ${\mathfrak S}_{n-1}\subset{\mathfrak S}_n$ (the standard notation is $\Res_{{\mathfrak S}_{n-1}}^{{\mathfrak S}_n}V_\lambda)$ is the sum of $s$ irreducible representation of ${\mathfrak S}_{n-1}$: $\Res_{{\mathfrak S}_{n-1}}^{{\mathfrak S}_n}V_\lambda=V_{\lambda\downarrow_1}\oplus\ldots\oplus V_{\lambda\downarrow_s}$. Hence,$$\dim V_\lambda=\dim V_{\lambda\downarrow_1}+\dots+\dim V_{\lambda\downarrow_s}.\eqno(3)$$

It is worth mentioning, although we do not need it in the current article, that for Young diagrams there exists an operation, in some sense ``dual" to removing removable cells, namely attaching attachable cells. A cell {\sl attachable} to the Young diagram is a cell {\sl not in} this diagram, such that after attaching this cell the Young diagram remains a Young diagram. It is easy to understand that the number of attachable cells for a Young diagram is one more than the number of removable cells. The partitions whose Young diagrams  are obtained from the Young diagram of $\lambda$ by attaching one attachable cell are denoted as $\lambda\hskip-3pt\uparrow^1,\lambda\hskip-3pt\uparrow^2,\dots \lambda\hskip-3pt\uparrow^{s+1}$. There is a well-known formula for the representation of ${\mathfrak S}_{n+1}$ induced by the representation $V_\lambda$ of ${\mathfrak S}_n$: $\Ind_{{\mathfrak S}_n}^{{\mathfrak S}^{n+1}}V_\lambda=V_{\lambda\uparrow^1}\oplus\dots\oplus V_{\lambda\uparrow^{s+1}}$. The similarity of this formula to the formula for $\Res_{{\mathfrak S}_{n-1}}$ given above, reflects {\it duality} between the operations Res and Ind (which can be expressed in terms of {\it adjoint functors}).

\section{Known results.}\label{known}

The results stated in three propositions below have been considered well known for a long time (see, for example, \cite{K2}). The first one is contained (and called ``Division theorem'') in the articles \cite{B1},\cite{B2} by A. Borodin. For the last two we prefer to give brief proofs. 

Let $\lambda=(\lambda_1,\dots,\lambda_N),\, \lambda_1\ge\dots\ge\lambda_N>0,\, \lambda_1+\dots+\lambda_N=n$ be a partition of $n$. We use the notations of Section 2: $\lambda'=(\lambda'_1,\dots,\lambda'_{N'})$ is a partition dual to $\lambda$, the number of removable cells in the Young diagram of $\lambda$ is $s$, and their coordinates are $(i_j,i'_j),\, j=1,\dots,s$, the partition whose Young diagram is obtained from Young diagram of $\lambda$ be removing the cell $(i_j,i'_j)$ is $\lambda_{\lambda\downarrow_j}$. In addition to this, we use the notation $\lambda'_0=\infty$. When necessary, we extend notations $n$ and $N$ to $n(\lambda)$ and $N(\lambda)$.

\begin{proposition}\label{P-recursion}$$P_\lambda(q)=\sum_{j=1}^s\left(q^{n-\lambda'_{i'_j}}-q^{n-1-\lambda'_{i'_j-1}}\right)P_{\lambda\downarrow_j}(q).\eqno(4)$$
\end{proposition}
For example, let $\lambda=(2,1)$. Then $\lambda'=\lambda=(2,1)$, $s=2$, the removable cells are $(1,2)$ and $(2,1),\ \lambda\hskip-3pt\downarrow_1=(1,1),\, \lambda\hskip-3pt\downarrow_2=(2)$. Obviously, $P_{(1,1)}(q)=1$, and $P_{(2)}=Q$. By (4), $P_{(2,1)}(q)=(q^{3-1}-q^{3-1-2})P_{(1,1)}(q)+(q^{3-2}-q^{3-1-\infty})P_{(2)}(q)=(q^2-1)\cdot1+qQ=Q(1+2q)$.

\begin{proposition}\label{degree}{\bf(a)} $$\deg P_\lambda(q)=\sum_{i<j}\lambda'_i\lambda'_j.\eqno(5)$$

{\bf(b)} The leading coefficient of the polynomial $P_\lambda(q)$ is equal to $\dim V_\lambda$.
\end{proposition}

{\it Proof.} Induction with respect to $n$. For $n=2$ the statement follows from the computations after Proposition \ref{P-recursion}. 

Put $D(\lambda)=\sum_{i<j}\lambda'_i\lambda'_j$. We want to prove that $\deg P_\lambda(q)=D(\lambda)$ and we assume known that $\deg P_{\lambda\downarrow_j}(q)=D(\lambda\hskip-3pt\downarrow_j)$ for every $j$ and that the leading coefficient of every $P_{\lambda\downarrow_j}(q)$ is $\dim V_{\lambda\downarrow_j}$. When we pass from $\lambda$ to $\lambda\hskip-3pt\downarrow_j$, all $\lambda'_i$ stay unchanged, except $\lambda'_{i'_j}$, which looses 1. Hence, $D(\lambda)-D(\lambda\hskip-3pt\downarrow_j)$ is the sum of all $\lambda'_i$, except $\lambda'_{i'_j}$, in other words, $D(\lambda)-D(\lambda\hskip-3pt\downarrow_j)=N-\lambda'_{i'_j}$. This shows that all the summands in the right hand part of (4) have the same degree $D(\lambda)$. Moreover, by the induction hypothesis, the leading coefficients of these summands are dimensions of $V_{\lambda\downarrow_j}$. Thus $P_\lambda$ has the degree $D(\lambda)$ and the leading coefficient $\sum_j\dim V_{\lambda\downarrow_j}$, which is,  by (3), $\dim V_\lambda$.

\begin{proposition}\label{e(lambda)}{\bf (a)} $e(\lambda)=n(\lambda)-N(\lambda)\strut$.

{\bf(b)} All the coefficients of the polynomial $Q^{-e(\lambda)}P_\lambda(q)$ (and hence of the polynomial $R_\lambda(q)$) are positive integers.
\end{proposition}

{\it Proof.} Again, the induction with respect to $n$. For $n=2$  the statement again follows from the computations after Proposition \ref{P-recursion}.

We assume that the statement is true for $\lambda\hskip-3pt\downarrow_j$. Obviously, $n(\lambda\hskip-3pt\downarrow_j)=n(\lambda)-1$ and $$N(\lambda\hskip-3pt\downarrow_j)=\left\{\begin{array} {ll} N(\lambda),&{\rm if}\ \lambda_N\ne1,\\ N(\lambda)-1&{\rm if}\ \lambda_N=1.\end{array}\right.$$Hence, $$\displaystyle{e(\lambda\hskip-3pt\downarrow_j)=\left\{\begin{array} {ll} n(\lambda)-N(\lambda)-1,&{\rm if}\ \lambda_N\ne1,\\ n(\lambda)-N(\lambda),&{\rm if}\ \lambda_N=1.\end{array}\right.}$$On the other hand, for the multiple in (4) we have:\vfil\eject
$$q^{N-\lambda'_{i'_j}}-q^{N-1-\lambda'_{i'_j-1}}=\left\{\begin{array} {ll}q^{N-1-\lambda'_{i'_j-1}}\cdot Q\left(1+q+\dots+q^{\lambda'_{i'_j}-\lambda'_{i'_j-1}}\right),&{\rm if}\ i'_j>1,\\ q^{N-\lambda'_{i'_j}},&{\rm if}\ i'_j=1\end{array}\right.$$\vskip-.1in $$\hskip1in=\left\{\begin{array} {ll} Q\cdot{\rm polynomial\ with\ positive \ integer\ coefficients},&{\rm if}\ i'_j>1,\\ \hskip.25in{\rm polynomial\ with\ positive \ integer\ coefficients,}&{\rm if}\ i'_j=1.\end{array}\right.$$ 

Notice that $\lambda_N=1$ and $i'_j=1$ are the same event: it means that the bottom cel in the right column of the Young diagram is removable, that is, the Young diagram has the shape as shown in Figure 2. Also, in this case, $j=s$.\vskip.2in

\centerline{
\beginpicture
\setcoordinatesystem units <1.2in,1.2in> point at 0 0
\linethickness=.6pt
\setplotsymbol({\bf.})
\plot .2 -1.1 .2 -1.4 0 -1.4 0 0 1.08 0 /
\setdashpattern <4pt,3.5pt>
\ellipticalarc axes ratio 2:1 -40 degrees from .2 -1.1 center at 1.7556 -1.1
\ellipticalarc axes ratio 2:1  -45 degrees from 1.08 0 center at -.56 0
\setsolid
\putrule from 0 -1.2 to .2 -1.2
\putrule from 0 -1 to .17 -1
\putrule from 0 -.8 to .26 -.8
\putrule from 0 -.6 to .46 -.6
\putrule from 0 -.4 to .76 -.4
\putrule from 0 -.2 to .98 -.2
\putrule from .2 0 to .2 -.9
\putrule from .4 0 to .4 -.66
\putrule from .6 0 to .6 -.52
\putrule from .8 0 to .8 -.36
\putrule from 1 0 to 1 -.16
\setplotsymbol({\tiny.})
\plot 0 -1.36 .04 -1.4 /
\plot 0 -1.32 .08 -1.4 /
\plot 0 -1.28 .12 -1.4 /
\plot 0 -1.24 .16 -1.4 /
\plot 0 -1.2 .2 -1.4 /
\plot .04 -1.2 .2 -1.36 /
\plot .08 -1.2 .2 -1.32 /
\plot .12 -1.2 .2 -1.28 /
\plot .16 -1.2 .2 -1.24 /
\put{Figure 2} at .1 -1.6
\endpicture
}\vskip.2in

Thus, in all cases, $P_\lambda(q)$ is $Q^{n(\lambda)-N(\lambda)}$ times a polynomial with positive integer coefficients. Since a polynomial with positive integer coefficients cannot be divisible by $Q$, this is what we had to prove.

\begin{corollary}$$\displaystyle{e(\lambda\hskip-3pt\downarrow_j)=\left\{\begin{array} {ll} e(\lambda)-1,&{\rm if}\ \lambda_N\ne1,\\ e(\lambda),&{\rm if}\ \lambda_N=1.\end{array}\right.}$$\end{corollary}

\section{Computation of $d(\lambda)$.}\label{dlambda}

We keep all the notation from Section \ref{known}; in addition to that, we put $\lambda'_{N'+1}=0$.

\begin{theorem} $$d(\lambda)={n(\lambda)\choose2}-{N(\lambda)\choose2}-\lambda'_1\lambda'_2-\lambda'_2\lambda'_3-\dots-\lambda'_{N'-1}\lambda'_{N'}\eqno(6)$$
\end{theorem}

{\bf Proof.} Until we have finished the proof, we denote the right hand side of (6) as $d'(\lambda)$. (Sometimes, we will append to the expression for $d'(\lambda)$ the zero term $-\lambda'_{N'}\lambda'_{N'+1}$.) Thus, we are proving that $d(\lambda)=d'(\lambda)$. 

We use the induction with respect to $n$. For $n=2$ the statement is true. Indeed, $P_{(2)}=Q,P_{(11)}=1,$ (see remark after Proposition \ref{P-recursion}), hence $d(2)=d(11)=0.$ Also, $(2)'=(11), (11)'=(2)$, hence

$$\begin{array} {c} d'(2)\displaystyle{={2\choose2}-{1\choose2}-1\cdot1=1-0-1=0=d(2);}\\d'(11)\displaystyle{={2\choose2}-{2\choose2}-2\cdot0=1-1-0=0=d(11).}\end{array}$$

By the induction hypothesis, (6) holds for all $\lambda\hskip-3pt\downarrow_j$. Let us compare $d(\lambda\hskip-3pt \downarrow_j)$ with $d'(\lambda)$. 

First, $n(\lambda\hskip-3pt\downarrow_j)=n(\lambda)-1$. Hence, $\displaystyle{n\choose2}$ in $d'(\lambda)$ becomes $\displaystyle{n-1\choose2}$, that is, looses $n-1$ in  $d(\lambda\hskip-3pt \downarrow_j)$. Second, $N(\lambda\hskip-3pt \downarrow_j)=\displaystyle \left\{\begin{array} {ll} N&{\rm if}\ i'_j\ne1,\\ N-1&{\rm if}\ i'_j=1.\end{array}\right.$ Hence, $-\displaystyle{N\choose2}$ becomes $-\displaystyle{N-1\choose2}$, hence, gains $N-1$, if $i'_j=1$. Third, when we pass from $\lambda$ to $\lambda\hskip-3pt_j$, all $\lambda'_i$'s remain the same, with the exception of $\lambda'_{i'_j}$, which becomes one less. Hence, $-\lambda'_{i'_j-1}\lambda'_{i'_j}-\lambda'_{i'_j}\lambda'_{i'_j+1}$ gains $\lambda'_{i'_j-1}+\lambda'_{i'_j+1}$ except the cases $i'_j=1$, when it gains just $\lambda'_2$, and $i'_j=N'$ when it gains just $\lambda'_{N'-1}$. Thus,$$\begin{array} {rll} d(\lambda\hskip-3pt\downarrow_1)&=d'(\lambda)-(n-1)+\lambda'_{i'_1-1},\\ d(\lambda\hskip-3pt\downarrow_j)&=d'(\lambda)-(n-1)+\lambda'_{i'_j-1}+\lambda'_{i'_j+1},&{\rm if}\ j>1\ {\rm and}\ i'_j\ne1,\\d(\lambda\hskip-3pt\downarrow_s)&=d'(\lambda)-(n-1)+(N'-1)+\lambda'_2,&{\rm if}\ i'_s=1.\end{array}\eqno(7)$$

Remind that $d(\lambda\hskip-3pt \downarrow_j)$ is the lowest exponent at $q$ in the polynomial $P_{\lambda\downarrow_j}$. According to Proposition 3.1, $P_\lambda$ is the sum of the polynomials $P_{\lambda\downarrow_j}$ multiplied, respectively, by $q^{n-\lambda'_{i'_j}}-q^{n-1-\lambda'_{i'_j-1}}$. This multiplication increases the lowest exponent at $q$ in $P_{\lambda\downarrow_j}$ by $n-1-\lambda'_{i'_j-1}$, if $i'_j\ne1$ and increases it by $n-\lambda'_1=n-N$, if $i'_j=1$. This increased exponent is $d'(\lambda)$, if $j=1$ and $d'(\lambda)+\lambda'_{i'_j+1}$ in all other cases. Thus, $P_\lambda$ is the sum of polynomials, one of which has the lowest exponent at $q$ equal to $d'(\lambda)$, while all the rest have this lowest exponent greater than $d'(\lambda)$. Hence the lowest exponent at $q$ in $P_\lambda$, that is $d(\lambda)$, is equal to $d'(\lambda)$, which completes the proof.

\section{Corollaries for $R_\lambda$.}\label{Rlambda}

Results of Sections \ref{known} and \ref{dlambda} may be used to calculate the degree of the polynomial $R_\lambda$:\smallskip

\begin{theorem}\label{degR} $$\deg R_\lambda=\sum_{i=1}^{N'-1}\lambda'_i\lambda'_{i+1}-\sum_{i=2}^{N'}{\lambda'_i+1\choose2}. \eqno(8)$$ \end{theorem}

{\bf Proof.} First, we need to transform the expression (4) for $\deg P_\lambda$: $$\begin{array} {rl} \deg P_\lambda&=\displaystyle{\sum_{i<j}\lambda'_i\lambda'_j=\frac{\left(\sum\nolimits_i\lambda'_i\right)^2-\sum\nolimits_i\left(\lambda'_i\right)^2}2}\\ &=\displaystyle\frac{\left(\sum\nolimits_i\lambda'_i\right)^2-\sum_i\lambda'_i\left(\lambda'_i-1\right)-\sum_i\lambda'_i}2\\ &=\displaystyle{\frac{n^2-n}2-\sum_{i=1}^{N'}\frac{\lambda'_i(\lambda'_i-1)}2={n\choose2}-\sum_{i-1}^{N'}{\lambda'_i\choose2}}.\end{array}$$Using this and expressions for $d(\lambda)$ and $e(\lambda)$ from Theorem 4.1 and Proposition 3.3(a), we obtain:

$$\begin{array} {rl} \deg R_\lambda&=\deg P_\lambda-d(\lambda-e(\lambda)\strut\\ &\displaystyle{=\left[{n\choose2}-\sum_{i=1}^{N'}{\lambda'_i\choose2}\right]-\left[{n\choose2}-{N\choose2}-\sum_{i-1}^{N'-1}\lambda'_i\lambda'_{i+1}\right]-(n-N)}\\ &=\displaystyle{\sum_{i=1}^{N'-1}\lambda'_i\lambda'_{i+1}+{N\choose2}+N-\sum\nolimits_{i=1}^{N'}{\lambda'_i\choose2}-n}\\&=\displaystyle{\sum_{i=1}^{N'-1}\lambda'_i\lambda'_{i+1}+{N\choose2}+N-\sum\nolimits_{i=1}^{N'}{\lambda'_i\choose2}-\sum\nolimits_{i=1}^{N'}\lambda'_i}\\&=\displaystyle{\sum_{i=1}^{N'-1}\lambda'_i\lambda'_{i+1}+{N+1\choose2}-\sum_{i=1}^{N'}{\lambda'_i+1\choose2}}\\&=\displaystyle{\sum_{i=1}^{N'-1}\lambda'_i\lambda'_{i+1}+{\lambda'_1+1\choose2}-\sum_{i=1}^{N'}{\lambda'_i+1\choose2}=\sum_{i=1}^{N'-1}\lambda'_i\lambda'_{i+1}-\sum_{i=2}^{N'}{\lambda'_i+1\choose2}} \end{array}$$

Next, we can derive a formula similar to (4) for $R_\lambda$:\smallskip

\begin{theorem}\label{R-recursion}  $$R_\lambda=\sum_{j=1}^sq^{\lambda'_{i'_j+1}}\left(1+q+\dots+q^{\alpha_{i'_j-1}}\right)R_{\lambda\downarrow_j}={\sum_{j=1}^s\frac{q^{\lambda'_{i'_j+1}}\left(1-q^{\alpha_{i'_j-1}+1}\right)}{1-q}}R_{\lambda\downarrow_j}\eqno(9)$$\emph{(if $j=1$, then $i'_j=N'$, and $\lambda'_{i'_j+1}=\lambda'_{N'+1}=0$; if $i'_j=1$, then the sum in parentheses in the first formula is 1)}.\end{theorem}

{\bf Example:} $\displaystyle{R_{2^\alpha1^\beta}=q^\alpha R_{2^\alpha1^{\beta-1}}+\left(1+q+\dots+q^\beta\right)R_{2^{\alpha-1}1^{\beta-1}}}.$\smallskip

{\bf Proof of Theorem 5.2.} Replace in formula (4) $P_\lambda$ and $P_{\lambda\downarrow_j}$ by, respectively, $q^{d(\lambda)}Q^{e(\lambda)}R_\lambda$ and $q^{d(\lambda\downarrow_j)}Q^{e(\lambda\downarrow_j)}R_{\lambda\downarrow_j}$, and then divide the equality by $q^{d(\lambda)}Q^{e(\lambda)}$. We get:

$$R_\lambda(q)=\sum_{j=1}^s\left(q^{n-\lambda'_{i'_j}}-q^{n-1-\lambda'_{i'_j-1}}\right)q^{d(\lambda\downarrow_j)-d(\lambda)}Q^{e(\lambda\downarrow_j)-e(\lambda)}R_{\lambda\downarrow_j}(q).$$

By Corollary in Section \ref{known}, $\displaystyle{Q^{e(\lambda\downarrow_j)-e(\lambda)}=\left\{\begin{array} {ll}Q^{-1}&{\rm if}\ i'_j\ne1,\\ 1&{\rm if}\ i'_j=1.\end{array}\right.}$The difference $d\left(\lambda\hskip-5pt\downarrow_j\right)-d(\lambda)$ is described in formulas (7) (in which $d'(\lambda)$ should be replaced by $d(\lambda)$). 

If $i'_j\ne1$, then$$q^{n-\lambda'_{i'_j}}-q^{n-1-\lambda'_{i'_j-1}}=q^{n-1-\lambda'_{i'_j-1}}\left(q^{\lambda'_{i'_j-1}-\lambda'_{i'_j}}+1\right)=q^{n-1-\lambda'_{i'_j-1}}\left(q^{\alpha'_{j-1}+1}-1\right).$$Hence $$\begin{array} {l} \left(q^{n-\lambda'_{i'_j}}-q^{n-1-\lambda'_{i'_j-1}}\right)q^{d(\lambda\downarrow_j)-d(\lambda)}Q^{e(\lambda\downarrow_j)-e(\lambda)}\\ \hskip.6in=q^{n-1-\lambda'_{i'_j-1}}\displaystyle\frac{q^{\alpha'_{j-1}+1}-1}Qq^{-(n-1)+\lambda'_{i'_j-1}+\lambda'_{i'_j+1}}=q^{\lambda'_{i'_j+1}}\left(1+q+\dots+q^{\alpha_{i'_j-1}}\right)\end{array}$$

If $i'_j=1$, then $j=s,\, \lambda'_{i'_j}=N'$, and $$\left(q^{n-\lambda'_{i'_j}}-q^{n-1-\lambda'_{i'_j-1}}\right)q^{d(\lambda\downarrow_j)-d(\lambda)}Q^{e(\lambda\downarrow_j)-e(\lambda)}=q^{n-N'}q^{-(n-1)+{n'-1}+\lambda'_2}=q^{\lambda'_2}\ \left(=q^{\lambda'_{i'_j+1}}\right).$$

This completes the proof.

\begin{corollary}\label{constant} The constant term of the polynomial $R_\lambda$ is $1$.\end{corollary}

This follows from (9): the only summand in this equality, which can have a non-zero constant term is that with $j=1$, and this constant term is that of $R_{\lambda\downarrow_s}$. Our statement follows by induction.

\section{Stabilization of $\bf R_\lambda$.}\label{stabilization}

For $\lambda=\dots3^{\alpha_3}2^{\alpha_2}1^{\alpha_1}$ put $\lambda^+=\dots3^{\alpha_3}2^{\alpha_2}1^{\alpha_1+1}$.

\begin{theorem}\label{plus1} The polynomials $R_{\lambda^+}$ and $R_\lambda$ have equal coefficients at $1,q,\dots,q^{\alpha_1}$. In other words, the polynomial $R_{\lambda^+}-R_\lambda$ is divisible by $q^{\alpha_1+1}$.\end{theorem}

 {\bf Proof.}  If $\alpha_1=0$, then we need to prove only that the polynomials $R_{\lambda^+}$ and $R_\lambda$ have equal constant terms, but we know this from Corollary \ref{constant}. Thus, we can assume that $\alpha_1>0$. In this case it is obvious that $\lambda^+\hskip-5pt\downarrow_j=(\lambda\hskip-3pt\downarrow_j)^+$, so we can use the notation $\lambda\hskip-3pt\downarrow_j^+$.
 
One more trivial case: if $\lambda=1^n$, then $\lambda^+=1^{n+1}$ and $R_\lambda=R_{\lambda^+}=1$, so our statements trivially holds. Thus, we can assume that $s>1$. 

We use the induction with respect to $n$ (for small values of $n$ our statement is obviously true). Thus, we can assume that the statement is true for all $\lambda\hskip-3pt\downarrow_j$'s. For all partitions $\lambda\hskip-3pt\downarrow_j$ with $j<s$ the number of ones is (at least) $\alpha_1$. For $\lambda\hskip-3pt\downarrow_s$, it is one less: $\alpha_1-1$. 

Now let us compare expressions (9) for $R_\lambda$ and $R_{\lambda^+}$. The number $s$ of summands will be the same.  The factors $q^{\lambda'_{i'_j+1}}\left(1+q+\dots+q^{\alpha_{i'_j-1}}\right)$ for the summands with $j<s$ will also be the same. As to $R_{\lambda\downarrow_j}$ and  $R_{\lambda\downarrow_j^+}$, their difference is divisible by $q^{\alpha_1+1}$. Thus, it remains to compare the terms with $j=1$. The polynomials $R_{\lambda\downarrow_1}$ and $R_{\lambda\downarrow_1^+}$, but this is compensated by the multiplication by $q^{\lambda'_{i'_j+1}}$, which is a positive power of $q$. This completes the proof of Theorem \ref{plus1}.\smallskip

Theorem \ref{plus1} shows that the coefficients at every fixed degree of $q$ in the sequence $$R_\lambda,R_{\lambda^+},R_{\lambda^{++}},R_{\lambda^{+++}},\dots$$ stabilize, and in the limit we get a power series. This  series does not depend on 1's in $\lambda$, so we can assume that $\lambda$ has zero $\alpha_1$. For the limit series we will use the notation $R_{\lambda1^\infty}$. We will add to removable cells $(i_1,i'_1),\dots,(i_s,i'_s)$ of the Young diagram of $\lambda$ an imaginary ``removable cell" $(i_{s+1},i'_{s+1})=(\infty,1)$. There arise $s+1$ ``partitions" $\lambda1^\infty\hskip-3pt\downarrow_j$, which are $\lambda\hskip-3pt\downarrow_1\hskip-3pt1^\infty,\dots,\lambda\hskip-3pt\downarrow_s\hskip-3pt1^\infty,$  and $\lambda1^\infty\hskip-3pt\downarrow_{s+1}=\lambda1^\infty$. (It is possible that the partition $\lambda\hskip-3pt\downarrow_s$ has a non-zero $\alpha_1$, but we still use the notation $\lambda\hskip-3pt\downarrow_s\hskip-3pt1^\infty$.)

If $\lambda=(N')^{\alpha_N'}(N'-1)^{\alpha_{N'-1}}\dots3^{\alpha_3}2^{\alpha_2}$, then $\lambda1^\infty=(N')^{\alpha_{N'}}(N'-1)^{\alpha_{N'-1}}\dots3^{\alpha_3}2^{\alpha_2}1^\infty$ and the ``partition" $(\lambda1^\infty)'$ of $\infty$ is $(\infty,\lambda'_2,\dots,\lambda'_{N'}),$ where $\lambda'_i=\alpha_i+\alpha_{i+1}+\dots+\alpha_{N'}$. Thus for $i>1$, $(\lambda1^\infty)'_i=\lambda'_i$. 

In the statement below, we use the following notation: $S_k=(1-q)(1-q^2)\dots(1-q^k)$; in particular, $S_0=1$.

\begin{theorem}\label{stable} For $\lambda=2^{\alpha_2}3^{\alpha_3}\dots (N')^{\alpha_{N'}}$, $$R_{\lambda1^\infty}=\frac1{(1-q)^{n-N}S_{\alpha_2}\dots S_{\alpha_{N'}}}\eqno(10)$$\emph{(the two expressions for $R_{\lambda1^\infty}$ are the same, since the only non-zero $\alpha_i$'s  (for $\lambda$) are $\alpha_{i'_j}$'s -- see formula (1) in Section \ref{partitions}).}\end{theorem}

{\bf Proof.} Again, the induction with respect to $N$. Thus, we assume that formula (10) holds for $R_{\lambda\downarrow_j}$. Let us look, how much the (lower) expression (10) changes when we replace $\lambda$ by $\lambda\hskip-3pt\downarrow_j$. First, $N$ stays unchanged, but $n$ becomes 1 less, so $n-N$ becomes 1 less. Second, $\alpha_{i'_j}$ becomes 1 less, so $S_{\alpha_{i'_j}}$ in the denominator looses (that is, the numerator gains) the factor $1-q^{i'_j}$. Third, $\alpha_{i'_j-1}$ (which could be zero) becomes 1 bigger, so the denominator gains the factor $1-q^{\alpha_{i_j-1}+1}$. Let us summarize:$$R_{\lambda\downarrow_j1^\infty}=\frac1{(1-q)^{n-N}S_{\alpha_{i'_1}}\dots S_{\alpha_{i'_s}}}\cdot\frac{(1-q)\left(1-q^{\alpha_{i'_j}}\right)}{1-q^{\alpha_{i_j-1}+1}}.\eqno(11)$$

Now let us apply the formula (9) to the ``partition" $\lambda1^\infty$ (that is, apply it to the partition $\lambda1^K$ with a very large $K$ and take the limit for $K\to\infty$). Since for $\lambda1^\infty$ the number $s$ becomes $s+1$, we get the following:$$R_{\lambda^\infty}=\sum_{j=1}^{s+1}\frac{q^{(\lambda1^\infty)'_{i'_j+1}}\left(1-q^{\alpha_{i'_j-1}+1}\right)}{1-q}R_{(\lambda1^\infty)\downarrow_j}.\eqno(12)$$

First, let us consider a summand in the last sum for a $j\le s$; we will take into account the equalities $(\lambda1^\infty)\hskip-3pt\downarrow_j=\lambda\hskip-3pt\downarrow_j\hskip-3pt1^\infty$ and $(\lambda1^\infty)'_i=\lambda'_i$ (see  above), and formula (11). We have:$$\begin{array} {rl}&\hskip.8in\displaystyle{\frac{q^{(\lambda1^\infty)'_{i'_j+1}}\left(1-q^{\alpha_{i'_j-1}+1}\right)}{1-q}R_{(\lambda1^\infty)\downarrow_j}} =\displaystyle{\frac{q^{\lambda'_{i'_j+1}}\left((1-q^{\alpha_{i'_j-1}+1}\right))}{1-q}R_{\lambda\downarrow_j1^\infty}}\\
&=\displaystyle\frac{q^{\lambda'_{i'_j+1}}\left(1-q^{\alpha_{i'_j-1}+1}\right)}{1-q}\cdot
\displaystyle{\frac{(1-q)\left(1-q^{\alpha_{i'_j}}\right)}{1-q^{\alpha_{i_j-1}+1}}}\cdot\displaystyle\frac1{(1-q)^{n-N}S_{\alpha_{i'_1}}\dots S_{\alpha_{i'_s}}}\\ &=\displaystyle\frac{q^{\lambda'_{i'_j+1}}\left(1-q^{\alpha_{i'_j}}\right)}{(1-q)^{n-N}S_{\alpha_{i'_1}}\dots S_{\alpha_{i'_s}}}=\displaystyle\frac{q^{\lambda'_{i'_j+1}}-q^{\lambda'_{i'_j+1}+\alpha_{i'_j}}}{(1-q)^{n-N}S_{\alpha_{i'_1}}\dots S_{\alpha_{i'_s}}}=\displaystyle\frac{q^{\lambda'_{i'_j+1}}-q^{\lambda'_{i'_j}}}{(1-q)^{n-N}S_{\alpha_{i'_1}}\dots S_{\alpha_{i'_s}}}\\ \end{array}$$(for the last equality, we used the relation $\lambda'_{i'_j}=\lambda'_{i'_j+1}+\alpha_{i'_j}$.

Now, let us turn to the summand in the formula (12), which corresponds to $j=s+1$. Since $i'_{s+1}=1$,  $\alpha_0=0$, and $R_{(\lambda1^\infty)\downarrow_{s+1}}=R_{\lambda1^\infty}$, this term is$$\frac{q^{\lambda'_{1+1}}(1-q)}{1-q}R_{\lambda1^\infty}=q^{\lambda'_2}R_{\lambda1^\infty}.$$ The final result  is:$$R_{\lambda1^\infty}=\frac{q^{\lambda'_{i'_1+1}}-q^{\lambda'_{i'_1}}+q^{\lambda'_{i'_2+1}}-q^{\lambda'_{i'_2}}+\dots+q^{\lambda'_{i'_s+1}}-q^{\lambda'_{i'_s}}}{(1-q)^{n-N}S_{\alpha_{i'_1}}\dots S_{\alpha_{i'_s}}}+q^{\lambda'_2}R_{\lambda1^\infty}.\eqno(13)$$

Notice that $\lambda'_{i'_j}=\lambda'_{i'_{j+1}+1}$. Indeed, $$\lambda'_{i'_{j+1}+1}=\alpha_{i'_{j+1}+1}+\dots+\alpha_{i'_j-1}+\underbrace{\alpha_{i'_j}+\dots+\alpha_{N'}}_{\lambda'_{i'_j}}$$and$$\alpha_{i'_{j+1}+1}=\dots=\alpha_{i'_j-1}=0.$$Therefore, in the numerator in the fraction in (13) everything, except the first and the last term, cancels. As to these two terms, $q^{\lambda'_{i'_1+1}}=q^0=1$, and $q^{\lambda'_{i'_s}}=q^{\lambda'_{i'_{s+1}+1}}=q^{\lambda'2}$. Thus, we have:$$R_{\lambda1^\infty}=\frac{1-q^{\lambda'_2}}{(1-q)^{n-N}S_{\alpha_{i'_1}}\dots S_{\alpha_{i'_s}}}+q^{\lambda'_2}R_{\lambda1^\infty}\Longrightarrow R_{\lambda1^\infty}=\frac1{(1-q)^{n-N}S_{\alpha_2}\dots S_{\alpha_{N'}}},$$which is what we had to prove.\smallskip

{\bf Examples.}\smallskip

$\bullet$ If $i_1>i_2>\dots>i_k>1$, then $R_{(i_1,i_2,\dots,i_k,1^\infty)}=\displaystyle\frac1{(1-q)^{i_1+i_2+\dots+i_k}}$.

$\bullet$ $R_{N^k1^\infty}=\displaystyle\frac1{(1-q)^{N(k-1)+1}(1-q^2)(1-q^3)\dots(1-q^k)}$

$\bullet$ $R_{N^2(N-1)^2\dots3^22^21^\infty}=\displaystyle\frac1{(1-q)^{N^2-1}(1-q^2)^{N-1}}$.

\section{Appendix. Table for $P_\lambda,\, n(\lambda)\le10$.}

The table below is presented in two columns on the next page and in one column in subsequent pages. In every line, the first entry is the notation for the partition $\lambda$, the second entry is the part $q^{d(\lambda)}Q^{e(\lambda)}$ of the polynomial $P_\lambda$ and the sequence in square brackets is the sequence of coefficients of the polynomial $R(\lambda)$ (starting with the constant term. For example, the line$$31\ q^21Q^2\ [1,3]$$ means $P_{31}=q^2Q^2(1+3q)$, and the line $$51^3\ q^{15}Q^4\ [1,5,15,35]$$ means $P_{51^3}=q^{15}Q^4(1+5q+15q^2+35q^3).$\vskip.2in

\noindent
\hangindent=-3in \hangafter=-26 
$\begin{array} {rll}
{\bf2}&Q&[1]\\
1^2&1&[1]\\
{\bf3}&qQ^2&[1]\\
21&Q&[1,2]\\
1^3&1&[1]\\
{\bf4}&q^3Q^3&[1]\\
31&q^2Q^2&[1,3]\\
2^2&qQ^2&[1,2]\\
21^2&Q&[1,2,3]\\
1^4&1&[1]\\
{\bf5}&q^6Q^4&[1]\\
41&q^5Q^3&[1,4]\\
32&q^3Q^3&[1,4,5]\\
31^2&q^3Q^2&[1,3,6]\\
2^21&qQ^2&[1,3,6,5]\\
21^3&Q&[1,2,3,4]\\
1^5&1&[1]\\
{\bf6}&q^{10}Q^5&[1]\\
51&q^9Q^4&[1,5]\\
42&q^7Q^4&[1,5,9]\\
41^2&q^7Q^3&[1,4,10]\\
3^2&q^6Q^4&[1,4,5]\\
321&q^4Q^3&[1,5,14,24,16]\\
31^3&q^4Q^2&[1,3,6,10]\\
2^3&q^3Q^3&[1,3,6,5]\\
2^21^2&qQ^2&[1,3,7,12,13,9]\\
21^4&Q&[1,2,3,4,5]\\
1^6&1&[1]\\
{\bf7}&q^{15}Q^6&[1]\\61&q^{14}Q^5&[1,6]\\
52&q^{12}Q^5&[1,6,14]\\
\end{array}$
\vskip-6.25in
\hangindent=2.6in \hangafter=-30
$\begin{array} {rll}
51^2&q^{12}Q^4&[1,5,15]\\
43&q^{10}Q^5&[1,6,14,14]\\
421&q^9Q^4&[1,6,20,43,35]\\
41^3&q^9Q^3&[1,4,10,20]\\
3^21&q^8Q^4&[1,5,15,28,21]\\
32^2&q^6Q^4&[1,5,15,28,35,21]\\
321^2&q^5Q^3&[1,5,15,34,58,62,35]\\
31^4&q^5Q^2&[1,3,6,10,15]\\
2^31&q^3Q^3&[1,4,10,20,28,28,14]\\
2^21^3&qQ^2&[1,3,7,13,21,24,22,14]\\
21^5&Q&[1,2,3,4,5,6]\\
1^7&1&[1]\\
{\bf8}&q^{21}Q^7&[1]\\
71&q^{20}Q^6&[1,7]\\
62&q^{18}Q^6&[1,7,20]\\
61^2&q^{18}Q^5&[1,6,21]\\
53&q^{16}Q^6&[1,7,20,28]\\
521&q^{15}Q^5&[1,7,27,69,64]\\
51^3&q^{15}Q^4&1,5,15,35]\\
4^2&q^{15}Q^6&[1,6,14,14]\\
431&q^{13}Q^5&[1,7,27,69,106,70]\\
42^2&q^{12}Q^5&[1,6,21,48,78,56]\\
421^2&q^{11}Q^4&[1,6,21,55,112,135,90]\\
41^4&q^{11}Q^3&[1,4,10,20,35]\\
3^22&q^{10}Q^5&1,6,21,48,78,84,42]\\
3^21^2&q^{10}Q^4&[1,5,16,39,73,90,56]\\
32^21&q^7Q^4&[1,6,21,55,112,183,218,174,70]\\
321^3&q^6Q^3&[1,5,15,35,69,113,135,123,64]\\
31^5&q^6Q^2&[1,3,6,10,15,21]\\
2^4&q^6Q^4&1,4,10,20,28,28,14]\\
2^31^2&q^3Q^3&[1,4,11,24,45,68,87,88,64,28]\\
2^21^4&qQ^2&[1,3,7,13,22,33,39,39,33,20]\\
21^6&Q&[1,2,3,4,5,6,7]\\
1^8&1&[1]\\
\end{array}$\vfil\eject
$$\begin{array} {rll}
{\bf9}&q^{28}Q^8&[1]\\
81&q^{27}Q^7&[1,8]\\
72&q^{25}Q^7&[1,8,27]\\
71^2&q^{25}Q^6&[1,7,28]\\
63&q^{23}Q^7&[1,8,27,48]\\
621&q^{22}Q^6&[1,8,35,103,105]\\
61^3&q^{22}Q^5&[1,6,21,56]\\
54&q^{21}Q^7&[1,8,27,48,42]\\
531&q^{20}Q^6&[1,8,35,103,195,162]\\
52^2&q^{19}Q^6&[1,7,28,75,147,120]\\
521^2&q^{18}Q^5&[1,7,28,83,194,254,189]\\
51^4&q^{18}Q^4&[1,5,15,35,70]\\
4^21&q^{19}Q^6&[1,7,28,75,120,84]\\
432&q^{16}Q^6&[1,8,35,103,222,329,366,168]\\
431^2&q^{16}Q^5&[1,7,28,83,194,344,387,216]\\
42^21&q^{14}Q^5&[1,7,28,83,194,371,513,477,216]\\
421^3&q^{13}Q^4&[1,6,21,56,125,237,312,313,189]\\
41^5&q^{13}Q^3&[1,4,10,20,35,56]\\
3^3&q^{15}Q^6&[1,6,21,48,78,84,42]\\
3^221&q^{12}Q^5&[1,7,28,83,194,371,561,633,474,168]\\
3^21^3&q^{12}Q^4&[1,5,16,40,85,152,208,213,120]\\
32^3&q^{10}Q^5&[1,6,21,56,117,198,273,288,216,84]\\
32^21^2&q^8Q^4&[1,6,22,61,141,277,472,672,793,720,453,162]\\
321^4&q^7Q^3&[1,5,15,35,70,125,196,243,253,212,105]\\
31^6&q^7Q^2&[1,3,6,10,15,21,28]\\
2^41&q^6Q^4&[1,5,15,35,70,117,165,195,180,120,42]\\
2^31^3&q^3Q^3&[1,4,11,25,49,86,131,178,212,218,180,117,48]\\
2^21^5&qQ^2&[1,3,7,13,22,34,49,58,61,57,46,27]\\
21^7&Q&[1,2,3,4,5,6,7,8]\\
1^9&1&[1]\\
{\bf10}&q^{36}Q^9&[1]\\
91&q^{35}Q^8&[1,9]\\
82&q^{33}Q^8&[1,9,35]\\
81^2&q^{33}Q^7&[1,8,36]\\
73&q^{31}Q^8&[1,9,35,75]\\
721&q^{30}Q^7&[1,9,44,136,160]\\
71^3&q^{30}Q^6&[1,7,28,84]\\
64&q^{29}Q^8&[1,9,35,75,90]\\
631&q^{28}Q^7&[1,9,44,146,325,315]\\
\end{array}$$\vfil\eject

$$\begin{array} {rll}
62^2&q^{27}Q^7&[1,8,36,110,250,225]\\
621^2&q^{26}Q^6&[1,8,36,119,312,434,350]\\
61^4&q^{26}Q^5&[1,6,21,56,126]\\
5^2&q^{28}Q^8&[1,8,27,48,42]\\
541&q^{26}Q^7&[1,9,44,146,325,447,288]\\
532&q^{24}Q^7&[1,9,44,146,360,654,828,450]\\
531^2&q^{24}Q^6&[1,8,36,119,312,641,836,567]\\
52^21&q^{22}Q^6&[1,8,36,119,312,676,1036,1067,525]\\
521^3&q^{21}Q^5&[1,7,28,84,209,445,626,672,448]\\
51^5&q^{21}Q^4&[1,5,15,35,70,126]\\
4^22&q^{23}Q^7&[1,8,36,110,250,404,486,252]\\
4^21^2&q^{23}Q^6&[1,7,29,90,222,419,507,300]\\
43^2&q^{21}Q^7&[1,8,36,110,250,432,533,492,210]\\
4321&q^{19}Q^6&[1,9,44,154,431,988,1877,2838,3217,2304,768]\\
431^3&q^{19}Q^5&[1,7,28,84,209,445,791,1077,1033,525]\\
42^3&q^{18}Q^6&[1,7,28,84,200,392,644,801,693,300]\\
42^21^2&q^{16}Q^5&[1,7,29,90,231,507,973,1540,2026,2047,1432,567]\\
421^4&q^{15}Q^4&[1,6,21,56,126,251,446,601,676,616,350]\\
41^6&q^{15}Q^3&[1,4,10,20,35,56,84]\\
3^31&q^{18}Q^6&[1,7,28,84,200,392,609,711,558,210]\\
3^22^2&q^{15}Q^6&[1,7,29,90,222,454,782,1130,1338,1221,774,252]\\
3^221^2&q^{14}Q^5&[1,7,29,90,231,507,973,1615,2281,2647,2347,1422,450]\\
3^21^4&q^{14}Q^4&[1,5,16,40,86,165,281,395,461,425,225]\\
32^31&q^{11}Q^5&[1,7,28,84,209,445,826,1352,1918,2323,2323,1803,993,288]\\
32^21^3&q^9Q^4&[1,6,22,62,147,307,572,962,1432,1897,2186,2121,1635,935,315]\\
321^5&q^8Q^3&[1,5,15,35,70,126,209,315,395,435,420,334,160]\\
31^7&q^8Q^2&[1,3,6,10,15,21,28,36]\\
2^5&q^{10}Q^5&[1,5,15,35,70,117,165,195,180,120,42]\\
2^41^2&q^6Q^4&[1,5,16,40,86,165,281,430,591,725,775,710,525,285,90]\\
2^31^4&q^3Q^3&[1,4,11,25,50,90,150,225,310,390,449,461,409,310,190,75]\\
2^21^6&qQ^2&[1,3,7,13,22,34,50,69,82,88,87,78,61,35]\\
21^8&Q&[1,2,3,4,5,6,7,8,9]\\
1^{10}&1&[1]
\end{array}
$$

\end{document}